\documentclass{amsart}

\usepackage{amsthm,amsfonts,amsmath,amssymb,latexsym,epsfig}
\usepackage{upref,amssymb,eucal}

\begin{document}
\title{$p$-jets of $p$-isogenies}
\author{Alexandru Buium}
\def \bW{{\mathbb W}}
\def \bD{{\mathbb D}}
\def \cI{\mathcal I}
\def \cU{\mathcal U}
\def \cF{\mathcal F}
\def \cK{\mathcal K}
\def \cD{\mathcal D}
\def \cP{\mathcal P}
\def \cV{\mathcal V}
\def \cN{\mathcal N}
\def \cG{\mathcal G}
\def \cB{\mathcal B}
\def \cJ{\mathcal J}
\def \cI{\mathcal I}
\def \cA{\mathcal A}
\def \h{\hat{\ }}
\def \hp{\hat{\ }}
\def \tS{\tilde{S}}
\def \tP{\tilde{P}}
\def \tA{\tilde{A}}
\def \tX{\tilde{X}}
\def \tT{\tilde{T}}
\def \tE{\tilde{E}}
\def \tV{\tilde{V}}
\def \tC{\tilde{C}}
\def \tU{\tilde{U}}
\def \tG{\tilde{G}}
\def \tu{\tilde{u}}
\def \tx{\tilde{x}}
\def \tL{\tilde{L}}
\def \tY{\tilde{Y}}
\def \d{\delta}
\def \bZ{{\mathbb Z}}
\def \bV{{\mathbb V}}
\def \bF{{\bf F}}
\def \bE{{\bf E}}
\def \bC{{\bf C}}
\def \bO{{\bf O}}
\def \bR{\bmathbb R}
\def \bA{{\mathbb A}}
\def \bB{{\bf B}}
\def \cO{\mathcal O}
\def \ra{\rightarrow}
\def \bX{{\bf X}}
\def \bH{{\bf H}}
\def \bS{{\bf S}}
\def \bF{{\mathbb F}}
\def \bN{{\bf N}}
\def \bK{{\bf K}}
\def \bE{{\bf E}}
\def \bB{{\bf B}}
\def \bQ{{\bf Q}}
\def \bd{{\bf d}}
\def \bY{{\bf Y}}
\def \bU{{\bf U}}
\def \bL{{\bf L}}
\def \bQ{{\mathbb Q}}
\def \bP{{\bf P}}
\def \bR{{\bf R}}
\def \bC{{\mathbb C}}
\def \bM{{\mathbb M}}
\def \bG{{\mathbb G}}
\def \bP{{\bf P}}

\newtheorem{THM}{{\!}}[section]
\newtheorem{THMX}{{\!}}
\renewcommand{\theTHMX}{}
\newtheorem{theorem}{Theorem}[section]
\newtheorem{corollary}[theorem]{Corollary}
\newtheorem{lemma}[theorem]{Lemma}
\newtheorem{proposition}[theorem]{Proposition}
\theoremstyle{definition}
\newtheorem{definition}[theorem]{Definition}
\theoremstyle{remark}
\newtheorem{remark}[theorem]{Remark}
\theoremstyle{example}
\newtheorem{example}[theorem]{Example}
\numberwithin{equation}{section}
\address{Department of Mathematics and Statistics \\ University of New Mexico \\ Albuquerque, NM 87131, USA}
\email{buium@math.unm.edu}
\maketitle

\begin{abstract}
 $p$-jets of finite flat maps of schemes are generally neither finite nor flat.
  This phenomenon can be seen
 already in the case of $p$-isogenies of group schemes. However, for $p$-divisible groups,
 this pathology tends to disappear mod $p$ ``in the limit". We illustrate this in the case of the
 $p$-divisible groups of $\bG_m$, of elliptic curves, and of formal groups of finite height;
 in the  case of ordinary elliptic curves the structure in the limit  depends on the Serre-Tate parameter.
\end{abstract}

\section{Introduction}

Let $p$ be an odd prime and let $R=\widehat{\bZ_p^{ur}}$ be the $p$-adic completion of the maximum unramified extension of the ring $\bZ_p$ of $p$-adic numbers. (Throughout the paper the symbol $\widehat{\ }$ means $p$-adic completion.) Let $k=R/pR$ be the residue field of $R$.
 Then for each integer $n \geq 0$  a functor $J^n$ was introduced in \cite{char}
 that attaches to any scheme of finite type
$X/R$ a (Noetherian) $p$-adic formal scheme $J^n(X)$ over $R$ called the {\it $p$-jet space} of $X$ of order $n$.
For each $X$ there are  morphisms $J^n(X)\ra J^{n-1}(X)$, $n \geq 1$,  functorial in $X$, and $J^0(X)=\widehat{X}$.
We refer to \cite{char,book} for an exposition of the theory  and for some of the applications of these spaces; see also \cite{laplace, borger} for a several prime version of the theory.

 The functors $J^n$ behave nicely on smooth schemes and \'{e}tale morphisms:
 in particular if $X$ is smooth over $R$ then $J^n(X)$ are $p$-adic
 completions of smooth schemes over $R$; and if $X \ra Y$ is an \'{e}tale morphism  then $J^n(X)\simeq J^n(Y) \times_Y X$. So, in particular, if $X \ra Y$ is  finite and \'{e}tale  then the map $J^n(X) \ra J^n(Y)$ is again finite  and \'{e}tale.
 However if $X \ra Y$ is, say,  finite and flat  then $J^n(X)\ra J^n(Y)$ is generally neither finite nor flat. This basic
  pathology can be seen, in its simplest form, for $X$ a finite
  flat group scheme of degree  $p^{\nu}$ over $Y=Spec\ R$, $\nu\geq 1$; or  for $p$-isogenies
  $X \ra Y$ (i.e. isogenies of degree $p^{\nu}$) between smooth group schemes over $R$.
So a typical context for these phenomena is that of $p$-divisible
 groups $(X_{\nu};\nu\geq 1)$.
 The present paper offers an analysis of $p$-jet spaces of such groups.
 The main moral of the story will be that although the $p$-jet spaces of
 the individual $X_{\nu}$'s are generally
 highly pathological  order tends to be restored ``in the limit".

 Let us recall/introduce some basic notation. For any
 scheme of finite type $X/R$
  the rings of global functions $\cO^n(X):=\cO(J^n(X))$ form an
   inductive system; its direct limit is denoted by $\cO^{\infty}(X)$.
We typically denote by an upper bar reduction mod $p$;
in particular we set
\begin{equation}
\label{genneral}
\overline{\cO^n(X)}  := \cO^n(X)\otimes_R k,\ \ \overline{\cO^{\infty}(X)}  := \cO^{\infty}(X)\otimes_R k.\end{equation}
The morphisms $\overline{\cO^n(X)}\ra \overline{\cO^{\infty}(X)}$ are generally not injective (although they are injective for $X/R$ smooth \cite{char}).
 It turns out that for non-smooth $X/R$ a special role is then played by the rings:
 \begin{equation}
 \label{gennerall}
 \widetilde{\cO^n(X)}  :=  Im(\overline{\cO^n(X)}\ra \overline{\cO^{\infty}(X)}).
 \end{equation}
We will be mainly interested in the case when $X=G[p^{\nu}]$ are the kernels
 of the multiplication by $p^{\nu}$ on  $G=\bG_m$ (the multiplicative group),
 $G=E$ (an elliptic curve), or $G=\cF$ (a formal group of finite height).
  We  illustrate in this Introduction   the case of ordinary elliptic curves and that
  of formal groups of finite height.
(The remaining cases are special cases of these cases.)

So let first
$E/R$ be an elliptic curve with ordinary reduction and Serre-Tate parameter $q=q(E)\in U_1=1+pR$; cf. \cite{Katz}.
Also set $U_m=1+p^mR=U_1^{p^{m-1}}$, $m \geq 1$.
Let $$\theta\in \bQ_p/\bZ_p=\bigcup_{\nu=1}^{\infty} p^{-\nu}\bZ/\bZ$$ and assume $\nu_0\geq 1$ is minimal with the property
that $\theta \in p^{-\nu_0}\bZ/\bZ$.
For all $\nu\geq \nu_0$ let $E[p^{\nu}]$ be the kernel of the multiplication by $p^{\nu}$ in $E$, viewed as an extension of $p^{-\nu}\bZ/\bZ$ by $\mu_{p^{\nu}}$, and let $E[p^{\nu}]_{\theta}$ be the connected component of $E[p^{\nu}]$ lying over $\theta$. Also let $n \geq 1$. The following is part of Theorem \ref{painnn} in the body of the paper:

\begin{theorem}
\label{sppace}
\

1) If $q\not\in U_{\nu_0+1}$ and $\theta \neq 0$ then
$\widetilde{\cO^n(E[p^{\nu}]_{\theta})}=0$.

 2) If $q\in U_{\nu_0+1}$ or $\theta = 0$ then
$$
\widetilde{\cO^n(E[p^{\nu}]_{\theta})}   \simeq  \frac{k[x,x',x'',...,x^{(n)}]}{(x^{p^{\nu}},(x')^{p^{\nu}},...,(x^{(n)})^{p^{\nu}})}.$$
\end{theorem}

Here $x,x',...,x^{(n)}$ are variables and the isomorphisms in the theorem are compatible, in the obvious sense, with varying $n$ and $\nu$.
The statement of the theorem above should be contrasted with the fact that, as we shall see, for all $n,\nu\geq 1$ the
 the rings
 $\overline{\cO^n(E[p^{\nu}]_0)}$ have  positive Krull dimension (actually they are polynomial rings in  $\min\{n,\nu\}$  variables over some explicit local Artin rings).

Similarly let $\cF$ be a formal group law (in one variable) of
finite height  and let $\cF[p^{\nu}]$ be the kernel of the
multiplication by $p^{\nu}$ viewed as a finite flat group scheme
over $R$. Let $n \geq 1$. The following is part of Theorem \ref{mor}
in the body of the paper:

\begin{theorem}
\label{screw}
$$
\widetilde{\cO^n(\cF[p^{\nu}])}   \simeq
\frac{k[x,x',x'',...,x^{(n)}]}{(x^{p^{\nu}},(x')^{p^{\nu}},...,(x^{(n)})^{p^{\nu}})}.$$
\end{theorem}

Again we have compatibility with $n,\nu$ varying; also for $n,\nu\geq 1$,
 $\overline{\cO^n(\cF[p^{\nu}])}$ have  positive Krull dimension (and, again, they are polynomial rings in $\min\{n,\nu\}$  variables over some explicit local Artin rings).

\begin{remark}
It would be interesting to have a generalization
of Theorems \ref{sppace} and \ref{screw} to the case of  arbitrary $p$-divisible groups.
\end{remark}

\begin{remark}
It is interesting to note the following phenomenon. Let
$X_0,...,X_{\nu}$ be closed subschemes of the affine line $\bA^1$
over $R$ which are, say, finite and flat over $R$, and let
$$X=\bigcup_{i=0}^{\nu}X_i$$ (scheme theoretic union inside $\bA^1$,
defined by the intersection of the defining ideals). Then, in
general,
$$J^n(X) \neq \bigcup_{i=0}^{\nu} J^n(X_i)$$
as closed subschemes of $J^n(\bA^1)$. An example is provided by the
case when $X_i=Spec\ R[\zeta_{p^i}]$ where $\zeta_{p^i}$ is a
$p^i$-th root of unity. In this case $J^n(X_0)=Spec\ R$ and, as we
shall see later in the paper,  $J^n(X_i)=\emptyset$ for $i\geq 1$
and $n \geq 1$; on the other hand  $X=\mu_{p^{\nu}}$ (kernel of
multiplication by $p^{\nu}$ on $\bG_m$ over $R$) and hence
$J^n(\mu_{p^{\nu}})$ has a non-reduced reduction mod $p$ by Theorem
\ref{screw}. It would be interesting to understand this phenomenon
more generally when, for instance, $X_i=Spec\ R[\alpha_i]$ with
$\alpha_i$ integers in a finite ramified extension of the fraction
field of $R$.
\end{remark}

A few words about  proofs and the structure of the paper. The
structure of the rings in Theorems \ref{sppace} and \ref{screw} is
deceptively simple. What actually happens is that the equations
defining the $p$-jet spaces (of arbitrary order) of finite $p$-group
schemes (even of the simplest ones such as $\mu_{p^{\nu}}$) are a
priori extremely complex. The trick to handle this complexity (and
to eventually reveal the simple formulae in the above theorems) is
to introduce certain filtrations on the relevant rings and ideals
such that the Fermat quotient operator controlling the structure of
$p$-jet spaces is, in its turn, controlled by these filtrations in a
way that allows induction arguments to work. So after recalling from
\cite{char, book} some of the basic concepts we shall be dealing
with we will introduce the above mentioned filtrations and then, in
subsequent sections, we will investigate the $p$-jets of the
divisible groups of  $\bG_m$, ordinary elliptic curves, formal
groups of finite height, and supersingular elliptic curves
respectively. The $\bG_m$ case will be used as a step in the
analysis of  all the other cases. Also, the section on formal groups
 will contain more general results about $p$-jet spaces of schemes
defined by iterates of certain rather general power series that do
not necessarily come from formal groups.

\bigskip

{\it Acknowledgment}. This material is based upon work supported by
the National Science Foundation under Grant No. 0852591 and by the
IHES, France. Any opinions, findings, and conclusions or
recommendations expressed in this material are those of the author
and do not necessarily reflect the views of the National Science
Foundation or the IHES.

\section{Review of some basic concepts \cite{char,book}}

Rings in this paper will always be assumed commutative with unity.
A $p$-derivation $\d:A\ra A$ on a ring $A$ is a set theoretic map satisfying
$$\begin{array}{rcl}
\d(x+x) & = & \d x + \d y + C_p(x,y)\\
\d (xy) & = & x^p \d y + y^p \d x +p \d x \d y,\end{array}$$
where $C_p$ is the polynomial:
$$C_p(X,Y)=p^{-1}(X^p+Y^p-(X+Y)^p) \in \bZ[X,Y].$$
If $\d$ is as above then  $\phi:A \ra A$, $\phi(x)=x^p+p \d x$, is a ring homomorphism.
Note that $\d(xy)=x^p\d y+\phi(y)\d x=y^p\d x+\phi(x)\d y$.
Also $\d$ and $\phi$ commute. If $A$ is $p$-torsion free then $\d$ is, of course,  uniquely determined by $\phi$; also
\begin{equation}
\label{multiterm}
\d(x_1+...+x_m)=\d x_1+...+\d x_m+C_p(x_1,...,x_m),\end{equation}
 where
$$C_p(X_1,...,X_m):=p^{-1}(\sum_{i=1}^m X_i^p-(\sum_{i=1}^m X_i)^p)\in \bZ[X_1,...,X_m].$$

Now the ring $R=\widehat{\bZ_p^{ur}}$ has a unique $p$-derivation
defined by $\d x=(\phi(x)-x^p)/p$ where $\phi:R \ra R$ is the unique ring automorphism lifting the $p$-power Frobenius on $R/pR$.
Let $x$ be a variable (or more generally an $N$-tuple of variables $x_1,...,x_N$.)
  We consider the {\it $\d$-polynomial ring} $R\{x\}=R[x,x',x'',...]$; this is the polynomial ring in variables $x,x',x'',...,x^{(n)},...$, where $x',x'',...$ are variables (or $N$-tuples of variables), equipped with the unique $p$-derivation $\d:R\{x\}\ra R\{x\}$ such that $\d x=x'$, $\d x'=x''$, etc.
   For $X$ a scheme of finite type over $R$ one defines the $p$-jet spaces $J^n(X)$, $n \geq 0$ \cite{char}. The latter are $p$-adic formal schemes over $R$ fitting into a projective system
   $$...\ra J^n(X)\ra J^{n-1}(X) \ra ...\ra J^0(X)=\widehat{X}.$$
    Note that $X \mapsto J^n(X)$ are  functors commuting with open immersions and more generally with \'{e}tale maps in the sense that if $X \ra Y$ is \'{e}tale then $J^n(X)\simeq J^n(Y)\times_Y X$ in the category of $p$-adic formal schemes.
   If $X=Spec\ R[x]/(f)$ for a tuple of variables $x$ and a tuple of polynomials $f$ and then
   $$J^n(X)=Spf\ R[x,x',...,x^{(n)}]\widehat{\ }/(f,\d f,...,\d^n f).$$
   In particular if $Y \ra X$ is a closed immersion so is $J^n(Y)\ra J^n(X)$ for all $n$.
   Moreover $J^n$ commutes with fiber products: $J^n(Y\times_X Z)\simeq J^n(Y)\times_{J^n(X)}J^n(Z)$.
The rings $\cO^n(X):=\cO(J^n(X))$ form an inductive system and their direct limit $\cO^{\infty}(X)$ has a natural $p$-derivation $\d$ induced by the $p$-derivation $\d$ on $R\{x\}$. The following universality property holds. Assume for simplicity $X$ is affine. Then any $R$-algebra homomorphism of $\cO(X)$
into a $p$-adically complete ring $B$ equipped with a $p$-derivation $\d$ is induced by a unique $R$-algebra homomorphism $\cO^{\infty}(X)\ra B$ that commutes with $\d$.

 \section{Filtrations}

In this section we will introduce and study some basic filtrations, especially on rings equipped with $p$-derivations. In the next section we will specialize to the case when the ring in question is the ring $R\{x\}$ of $\d$-polynomials.

Let $A$ be a ring which for simplicity we assume $p$-torsion free and let $I$ be an ideal in $A$. For any integer $\nu \geq 0$ we denote by $I^{[p^{\nu}]}$ the ideal generated by all the elements of the form $p^if^{p^j}$ with $f\in I$, $i,j\geq0$, $i+j=\nu$. In particular $I^{[1]}=I$.

\begin{lemma}
\label{lem1}
\

1) $I^{[p^{\nu+1}]}\subset I^{[p^{\nu}]}.$

2) If $f \in I^{[p^{\nu}]}$ then $pf\in I^{[p^{\nu+1}]}$.

3) If $f \in I^{[p^{\nu}]}$ then $f^p\in I^{[p^{\nu+1}]}$.

4) If $I$ is generated by a family $\{f_s;s \in S\}$ then $I^{[p^{\nu}]}$ is generated by the family $\{p^i f_s^{p^j};s\in S, i,j\geq 0, i+j=\nu\}$.
\end{lemma}

{\it Proof}.
Assertions 1 and 2 are clear. For assertion 3 if $f=\sum_{s=1}^N p^{i_s}f_s^{p^{j_s}}g_s$ with $f_s\in I$, $g_s\in A$, and $i_s+j_s=\nu$ for all $s$ then
\begin{equation}
\label{maisus}f^p\in\sum_{s=1}^N (p^{i_s}f_s^{p^{j_s}}g_s)^p+
pC_p(p^{i_1}f_1^{p^{j_1}}g_1,...,p^{i_N}f_N^{p^{j_N}}g_N)
\subset I^{[p^{\nu+1}]}.\end{equation}
To check assertion 4 it is sufficient to prove that if $g\in I$ then $g^{p^{t}}$ is in the ideal
generated by the family $\{p^i f_s^{p^j};s\in S, i,j\geq 0, i+j=t\}$. One proves this by induction
 on $t \geq 0$. The case $t=0$ is clear. Now if the statement is true for $t=\nu$ and we set $f=g^{p^{\nu}}$ then we are done by equation (\ref{maisus}).
\qed

\bigskip

In what follows we assume we are given a $p$-derivation $\d:A \ra A$.

\begin{lemma}
\label{lem2}
Assume $\d (I)\subset I$. Then:

1) $\d(I^{[p^{\nu}]})\subset I^{[p^{\nu}]}$;

2) $\phi(I^{[p^{\nu}]})\subset I^{[p^{\nu+1}]}$.
\end{lemma}

{\it Proof}.
First,  for any $f \in A$, we have the following computation:
\begin{equation}
\label{deltafp}\begin{array}{rcl}
\d(f^{p^{\nu}}) & = & \frac{1}{p}(\phi(f^{p^{\nu}})-f^{p^{\nu+1}})\\
\  & = & \frac{1}{p}((f^p+p\d f)^{p^{\nu}}-f^{p^{\nu+1}})\\
\  & = & p^{\nu} f^{p(p^{\nu}-1)}\d f+p^{\nu+1}(\d f)^2 P(f,\d f)
\end{array}\end{equation}
where $P$ is a polynomial with $\bZ$-coefficients; indeed this is because
 for $2 \leq m \leq p^{\nu}$ we have
$$p^{\nu+1}\ |\ \left( \begin{array}{c} p^{\nu}\\ m\end{array} \right) p^{m-1}.$$
In particular if $f \in I$ then  $\d(f^{p^{\nu}})\in p^{\nu}I \subset I^{[p^{\nu}]}$.

 Let's prove assertion 1. In view of the equation (\ref{multiterm}) it is enough to note that for $i+j=\nu$, $f\in I$, $g\in A$ we have $\d(p^i f^{p^j}g)\in I^{[p^{\nu}]}$. Now
$$\begin{array}{rcl}
\d(p^i f^{p^j}g) & = & \d(p^i) f^{p^{j+1}}g^p + p^i \d(f^{p^j}g)\\
\  & = & \d(p^i)f^{p^{j+1}}g^p+p^i (\d(f^{p^j})g^p+f^{p^{j+1}}\d g +p (\d(f^{p^j}))(\d g))\\
\  & \in & I^{[p^{\nu}]},
\end{array}$$
by Lemma \ref{lem1} and because $\d(p^i)$ is either $0$ or in $p^{i-1}A$ according as $i=0$ or $i\geq 1$.

To prove assertion 2 note that if $f \in I^{[p^{\nu}]}$ then
$\phi(f)=f^p+p\d f\in I^{[p^{\nu+1}]}$
by Lemma \ref{lem1}.
\qed

\bigskip

In what follows we assume we are given, in addition,  a filtration on $A$
$$A^0\subset A^1 \subset A^2 \subset ...\subset A^n \subset ... \subset A,$$
by subrings $A^n$ such that $\d A^n \subset A^{n+1}$ for all $n \geq 0$. Then we define a new filtration
by subrings
\begin{equation}
\label{AAA}A^{\{0\}}\subset A^{\{1\}} \subset A^{\{2\}} \subset ...\subset A^{\{n\}} \subset ... \subset A,\end{equation}
$$A^{\{n\}}:=\sum_{s=0}^{\infty} p^s A^{n+s}=A^n+pA^{n+1}+p^2A^{n+2}+...$$

\begin{lemma}
\label{lem3}
\

1) $pA^{\{n+1\}}\subset A^{\{n\}}$;

2) $\d(A^{\{n\}})\subset A^{\{n+1\}}$;

3) $\phi(A^{\{n\}})\subset A^{\{n\}}$.
\end{lemma}

{\it Proof}. A trivial exercise. \qed

\bigskip

We will also need the following general:

\begin{lemma}
\label{binomial}
For any $f,g\in A^0$ we have the following equality of ideals in $A^n$:
$$(f-g,\d(f-g),...,\d^n(f-g))=(f-g,\d f -\d g, ...,\d^n f-\d^n g).$$
\end{lemma}

{\it Proof}.
Induction on $n$. The induction step follows from the congruence
$$\begin{array}{rclll}
\d(\d^{n-1} f -\d^{n-1} g)- (\d^n f -\d^n g) & = & C_p(\d^{n-1}f, -\d^{n-1}g) & \  &\ \\
\  & \equiv & C_p(\d^{n-1}g, -\d^{n-1}g) & \text{mod} & \d^{n-1}f-\d^{n-1}g\\
\  & = & 0. & \  & \  \end{array}$$
\qed

\section{$p$-jets of $p$-isogenies of $\bG_m$}

In this section we specialize the discussion of the previous section to the case when
\begin{equation}
\label{oz}
A=R\{x\},\ A^n=R[x,x',...,x^{(n)}],\ I=(x',x'',...).\end{equation}
So in this case, explicitly,
$$A^{\{n\}}=R[x,x',...,x^{(n)}]+pR[x,x',...,x^{(n+1)}]+p^2R[x,x',...,x^{(n+2)}]+...,$$
while $I^{[p^{\nu}]}$ is  the ideal of $R\{x\}$ generated by all $\d$-polynomials of the form
$p^i(x^{(s)})^{p^j}$,
with $s\geq 1$, $i,j\geq 0$, $i+j=\nu$; cf. assertion 4 in Lemma \ref{lem1}.

Here is the main result of this section:

\begin{lemma}
\label{coo}
For all integers $n,\nu\geq 1$ we have
$$\d^n(x^{p^{\nu}})\in
\begin{cases}
p^{\nu-n+1}x^{p^n(p^{\nu}-1)}\phi^{n-1}(x')+(p^{\nu-n+2}A^{\{1\}})\cap I^{[p^{\nu}]},
 & \text{ if $n \leq \nu+1$} \\
x^{p^n(p^{\nu}-1)}\phi^{\nu}(x^{(n-\nu)})+A^{\{n-\nu-1\}}\cap I^{[p^{\nu}]},
        & \text{ if $n \geq \nu+2$.}
\end{cases}
$$
\end{lemma}

{\it Proof}. First note that for $n \leq \nu+1$ we have
$\phi^{n-1}(x')  \in  I^{[p^{n-1}]}$ by Lemma \ref{lem2} and hence
\begin{equation}
\label{b1}
p^{\nu-n+1}x^{p^n(p^{\nu}-1)}\phi^{n-1}(x') \in I^{[p^{\nu}]},\end{equation}
by Lemma \ref{lem1}. Similarly, for $n \geq \nu+2$, we have
\begin{equation}
\label{b2}
x^{p^n(p^{\nu}-1)}\phi^{\nu}(x^{(n-\nu)})\in I^{[p^{\nu}]},\end{equation}
by Lemma \ref{lem2}. We also claim that
\begin{equation}
\label{b3}\d^n(x^{p^{\nu}})\in I^{[p^{\nu}]}.\end{equation}
To check (\ref{b3}) it is enough, by Lemma \ref{lem2}, to check that $\d(x^{p^{\nu}})\in I^{[p^{\nu}]}$; this however follows from equation (\ref{deltafp}). In view of (\ref{b1}), (\ref{b2}), (\ref{b3}), in order to prove our theorem it is enough to prove that
$$\d^n(x^{p^{\nu}})\in
\begin{cases}
p^{\nu-n+1}x^{p^n(p^{\nu}-1)}\phi^{n-1}(x')+p^{\nu-n+2}A^{\{1\}},
 & \text{ if $n \leq \nu+1$} \\
x^{p^n(p^{\nu}-1)}\phi^{\nu}(x^{(n-\nu)})+A^{\{n-\nu-1\}},
        & \text{ if $n \geq \nu+2$.}
\end{cases}
$$
We fix $\nu$ and proceed by induction on $n\geq 1$. For $n=1$ we are done by (\ref{deltafp}).
Next assume the Lemma is true for $n$ and we prove it for $n+1$.

Assume first $n \leq \nu +1$. By Lemma \ref{lem2}  $\phi^{n-1}(x')\in A^{\{1\}}$. So we have
$$\begin{array}{rcl}
\d^{n+1}(x^{p^{\nu}}) & \in & \d(p^{\nu-n+1}x^{p^n(p^{\nu}-1)}\phi^{n-1}(x'))+\d(p^{\nu-n+2}A^{\{1\}})\\
\  & \  & + C_p(p^{\nu-n+1}A^{\{1\}},p^{\nu-n+2}A^{\{1\}}).\end{array}
$$
Clearly the last term in the last equation is in $p^{\nu-n+1}A^{\{1\}}$. Also, by Lemma \ref{lem3}:
 $$\begin{array}{rcl}
\d(p^{\nu-n+2}A^{\{1\}}) & \subset & \d(p^{\nu-n+2})A^{\{1\}}+p^{\nu-n+2}\d(A^{\{1\}})\\
\  & \subset & p^{\nu-n+1}A^{\{1\}}+p^{\nu-n+2}A^{\{2\}}\\
\  & \subset & p^{\nu-n+1}A^{\{1\}}.
\end{array}$$

Now for $n \leq \nu$ we have
$$\d(p^{n-\nu+1})=p^{n-\nu}-p^{p(n-\nu+1)-1}\in p^{n-\nu}-p^{n-\nu+1}\bZ$$
hence:
$$\begin{array}{rcl}
\d(p^{\nu-n+1}x^{p^n(p^{\nu}-1)}\phi^{n-1}(x')) & = & \d(p^{\nu-n+1})x^{p^{n+1}(p^{\nu}-1)}
(\phi^{n-1}(x'))^p\\
\  & \  & +p^{\nu-n+1}\d(x^{p^n(p^{\nu}-1)}\phi^{n-1}(x'))\\
\  & \in & \ p^{\nu-n} x^{p^{n+1}(p^{\nu}-1)} \phi^{n-1}((x')^p)+p^{\nu-n+1}A^{\{1\}}\\
\  & \  & +p^{\nu-n+1}\d(x^{p^n(p^{\nu}-1)})(\phi^{n-1}(x'))^p\\
\  & \ & +p^{\nu-n+1} \phi(x^{p^n(p^{\nu}-1)})\d(\phi^{n-1}(x'))\\
\  & \subset & \ p^{\nu-n} x^{p^{n+1}(p^{\nu}-1)} \phi^{n-1}((x')^p)\\
\  & \  & +p^{\nu-n+1}A^{\{1\}}\\
\  & \ & +p^{\nu-n} (x^{p^{n+1}(p^{\nu}-1)}+pA^1)(\phi^{n-1}(px''))\\
\  & \subset & \ p^{\nu-n} x^{p^{n+1}(p^{\nu}-1)} \phi^{n}(x') +p^{\nu-n+1}A^{\{1\}}\end{array}$$
because $\d\circ \phi^{n-1}=\phi^{n-1}\circ \d$, $(x')^p+px''=\phi(x')$, and
$$p^{\nu-n}\cdot p  A^1 \cdot(\phi^{n-1}(px''))\subset p^{\nu-n+1} \cdot A^1 \cdot
 pA^{\{2\}}\subset p^{\nu-n+1} A^{\{1\}}.$$
So for $n\leq \nu$ we get
$$\d^{n+1}(x^{p^{\nu}})=  p^{\nu-n} x^{p^{n+1}(p^{\nu}-1)} \phi^{n}(x')+p^{\nu-n+1}A^{\{1\}},$$
which ends the induction step in  case $n \leq \nu$.

For $n=\nu+1$  we get
$$\begin{array}{rcl}\d(p^{\nu-n+1}x^{p^n(p^{\nu}-1)}\phi^{n-1}(x')) & = & \d(x^{p^n(p^{\nu}-1)}\phi^{n-1}(x'))\\
\  & =  & \d(x^{p^n(p^{\nu}-1)})(\phi^{n-1}(x'))^p\\
\  & \ & + \phi(x^{p^n(p^{\nu}-1)})\d(\phi^{n-1}(x'))\\
\  & \in & A^{\{1\}} + (x^{p^{n+1}(p^{\nu}-1)}+p A^1)(\phi^{n-1}(x''))\\
\  & = &  x^{p^{n+1}(p^{\nu}-1)} \phi^{n-1}(x'')+A^{\{1\}},\end{array}$$
by Lemma \ref{lem3}.
Hence
$$\d^{n+1}(x^{p^{\nu}})=x^{p^{n+1}(p^{\nu}-1)}\phi^{\nu}(x^{(n+1-\nu)})+A^{\{n-\nu\}},$$
which ends the induction step in case $n=\nu+1$.

Assume now $n \geq \nu+2$; then, by Lemma \ref{lem3},
$$
\begin{array}{rcl}
\d^{n+1}(x^{p^{\nu}}) & = & \d(x^{p^n(p^{\nu}-1)}\phi^{\nu}(x^{(n-\nu)}))+\d(A^{\{n-\nu-1\}})\\
\  & \  & + C_p(A^{\{n-\nu\}},A^{\{n-\nu-1\}})\\
\  & \in & \d(x^{p^n(p^{\nu}-1)}\phi^{\nu}(x^{(n-\nu)}))+A^{\{n-\nu\}}\\
 \ & = & x^{p^{n+1}(p^{\nu}-1)} \d(\phi^{\nu}(x^{(n-\nu)}))\\
\  & \  & + \phi^{\nu+1}(x^{(n-\nu)}) \d(x^{p^n(p^{\nu}-1)})+ A^{\{n-\nu\}}\\
\ & = & x^{p^{n+1}(p^{\nu}-1)} \phi^{\nu}(x^{(n+1-\nu)})\\
\  & \  & +  A^{\{n-\nu\}}\cdot pA^1 +A^{\{n-\nu\}}=A^{\{n-\nu\}}.
\end{array}
$$
This
 ends the induction step in case $n \geq \nu+2$.
\qed

\begin{corollary}
\label{dr}
$$\d^n(x^{p^{\nu}})\in \begin{cases}
A^{\{0\}}\cap I^{[p^{\nu}]} & \text{if $n \leq \nu$},\\
A^{\{n-\nu\}}\cap I^{[p^{\nu}]} & \text{if $n \geq \nu+1$}.
\end{cases}$$
\end{corollary}

For any element $f\in A^n$ we denote by $\overline{f}\in \overline{A}^n:=\overline{A^n}=k[x,x',...,x^{(n)}]$ the class of $f$ mod $p$. In particular one can consider the reduction mod $p$ of the ideal $I$, $\overline{I}=(x',x'',...)\subset \overline{A}=k[x,x',x'',...]$. Then
$$\overline{I^{[p^{\nu}]}}=\overline{I}^{[p^{\nu}]}=((x')^{p^{\nu}},(x'')^{p^{\nu}},...).$$
Moreover clearly
$$\overline{A}^n\cap \overline{I}^{[p^{\nu}]}=((x')^{p^{\nu}},...,(x^{(n)})^{p^{\nu}})\subset \overline{A}^n.$$
Also for any $f\in A^{\{n\}}\cap I^{[p^{\nu}]}$ we have $\overline{f} \in \overline{A}^n\cap\overline{I}^{[p^{\nu}]}$.

\begin{corollary} \label{modpp}
The element $\overline{\d^n(x^{p^{\nu}})}\in \overline{A}^n$
satisfies
$$\overline{\d^n(x^{p^{\nu}})}=
\begin{cases}
0, & \text{ if $1 \leq n\leq \nu$,}\\
x^{p^{\nu+1}(p^{\nu}-1)} (x')^{p^{\nu}}, & \text{ if $n=\nu+1$,}\\
x^{p^n(p^{\nu}-1)}(x^{(n-\nu)})^{p^{\nu}}+\text{(element in $\overline{A}^{n-\nu-1}\cap\overline{I}^{[p^{\nu}]}$)}, & \text{  if $n \geq \nu+2$.}\end{cases}$$
\end{corollary}

\begin{example}
The ``smallest" interesting case is $\nu=1$, $n=3$,
$$\overline{\d^3(x^{p})}=x^{p^3(p-1)}(x'')^{p}-\frac{1}{2} x^{p^3(p-2)}(x')^{2p}.$$
\end{example}

The following will also be useful later.

\begin{lemma}
\label{45}
\

1) $\d^n(x^p)\in R[x^p,x',...,x^{(n)}]$ for $n \geq 0$.

2) $\d^n(x_1x_2)\in R[x_1^p, x_2^p, x_1', x_2',...,x_1^{(n)}, x_2^{(n)}]$.
\end{lemma}

{\it Proof}.
Trivial induction on $n$.
\qed

\bigskip

 The natural context for Lemma \ref{coo} is as follows. Recall the
multiplicative group scheme
$$\bG_m=Spec\ R[x,x^{-1}]$$
over $R$ and the $p$-isogeny
$$[p^{\nu}]_{\bG_m}:\bG_m\ra \bG_m$$
given at the level of rings by $x \mapsto x^{p^{\nu}}$. The $p$-jet spaces of $\bG_m$ are given by
$$J^n(\bG_m)=Spf\ R[x,x^{-1},x',...,x^{(n)}]\h.$$
Consider the multiplication by $p^{\nu}$ map
$$[p^{\nu}]_{J^n(\bG_m)}:J^n(\bG_m)\ra J^n(\bG_m).$$
By the universality property of $J^n$ the map $[p^{\nu}]_{J^n(\bG_m)}$ coincides with the map $J^n([p^{\nu}]_{\bG_m})$
induced by $[p^{\nu}]_{\bG_m}$.
So the map $[p^{\nu}]_{J^n(\bG_m)}=J^n([p^{\nu}]_{\bG_m})$ is given at the level of rings by
\begin{equation}
\label{maap}
[p^{\nu}]^*:\cO^n(\bG_m)=R[x,x^{-1},x',...,x^{(n)}]\h\ra
R[x,x^{-1},x',...,x^{(n)}]\h,\end{equation}
$$
\ x\mapsto x^{p^{\nu}},\  x'\mapsto \d(x^{p^{\nu}}),\ ...,\ x^{(n)}\mapsto \d^n(x^{p^{\nu}}).$$
Lemma \ref{coo} can be seen then as a description of the map $[p^{\nu}]_{J^n(\bG_m)}=J^n([p^{\nu}]_{\bG_m})$.
Note that the map (\ref{maap}) is neither finite nor flat; indeed its reduction mod $p$
\begin{equation}
\label{barmaap} \overline{[p^{\nu}]^*}:\overline{\cO^n(\bG_m)}=
k[x,x^{-1},x',...,x^{(n)}]\ra
k[x,x^{-1},x',...,x^{(n)}]\end{equation}
 is neither finite nor flat because it is non-injective. See the Corollary \ref{aftert} below and also Proposition
\ref{abv} for more on this. Nevertheless note that by
 Corollary \ref{modpp} we have induced homomorphisms
$$k[x,x^{-1},x^{(\nu+1)},...,x^{(n)}]\simeq \frac{k[x,x^{-1},x',...,x^{(n)}]}{(x',...,x^{(\nu)})}\ra
k[x,x^{-1},x',...,x^{(n-\nu)}],$$ for $n \geq \nu+1$, and these
homomorphisms are finite algebra maps. By a dimension argument these
homomorphisms must be injective. Since finite maps of non-singular
varieties are automatically flat (\cite{E}, p. 461) we have:

\begin{corollary}
\label{aftert} For $n \geq \nu+1$ there are induced injective finite
flat maps
$$\overline{[p^{\nu}]^*}:\overline{\cO^n(\bG_m)}/(x',...,x^{(\nu)})\ra
\overline{\cO^{n-\nu}(\bG_m)}.$$
\end{corollary}

\section{$p$-jets of  $\mu_{p^{\nu}}$}

Start again with  the multiplicative group $\bG_m$ and the isogeny
$[p^{\nu}]_{\bG_m}:\bG_m \ra \bG_m$. Let
$$\mu_{p^{\nu}}=\bG_m[p^{\nu}]:=Ker([p^{\nu}]_{\bG_m})=Spec\ R[x,x^{-1}]/(x^{p^{\nu}}-1)=Spec\ R[x]/(x^{p^{\nu}}-1)$$
be the kernel of $[p^{\nu}]_{\bG_m}$. More generally, for any $a\in R^{\times}$, one may consider the finite flat scheme
$$\mu_{p^{\nu}}^a:=Spec\ R[x]/(x^{p^{\nu}}-a).$$
Its functor of points is given by
 $\mu_{p^{\nu}}^a(S)=\{s\in S; s^{p^{\nu}}=a\}$ for any $R$-algebra $S$.
Then $\mu_{p^{\nu}}^a$ has a natural structure of $\mu_{p^{\nu}}$-torsor.
More generally, for any $a,b \in R^{\times}$ we have  a natural morphism
$$\mu_{p^{\nu}}^a\times \mu_{p^{\nu}}^b \stackrel{can}{\longrightarrow} \mu_{p^{\nu}}^{ab}$$
given on $S$-points by $(s,t)\mapsto st$. We also have a natural isomorphism
$$\mu_{p^{\nu}}^{a^{p^{\nu}}b}\stackrel{can}{\longrightarrow} \mu_{p^{\nu}}^b$$
given on points by $s\mapsto s/a$.
Recall that we have set $U_m=1+p^{m}R=U_1^{p^{m-1}}$. So if $a\in U_{\nu+1}$ then $a=b^{p^{\nu}}$, $b \in U_1$, so division by $b$ gives an isomorphism
$$\mu_{p^{\nu}}^a\simeq \mu_{p^{\nu}}.$$

Finally the system $(\mu_{p^{\nu}};\nu \geq 1)$ is a $p$-divisible group with embeddings $\mu_{p^{\nu}} \subset  \mu_{p^{\nu+1}}$ given by the inclusions on points. More generally for any $a \in R^{\times}$ and any $\nu_0 \geq 1$ the schemes
$(\mu_{p^{\nu}}^{a^{p^{\nu-\nu_0}}};\nu \geq \nu_0)$,  form an inductive system  with embeddings
\begin{equation}
\label{omm}
\mu^a_{p^{\nu_0}}\subset \mu^{a^p}_{p^{\nu_0+1}} \subset ... \subset
\mu_{p^{\nu}}^{a^{p^{\nu-\nu_0}}}\subset \mu_{p^{\nu+1}}^{a^{p^{\nu+1-\nu_0}}}\subset ...\end{equation}
 given by the inclusions on points.
 Note that if $a\in U_{\nu_0+1}$ then we can write $a=b^{p^{\nu_0}}$ and hence division by $b$ gives an isomorphism between the inductive system
 (\ref{omm}) and the inductive system
 \begin{equation}
\label{ommm}
\mu_{p^{\nu_0}}\subset \mu_{p^{\nu_0+1}} \subset ... \subset
\mu_{p^{\nu}}\subset \mu_{p^{\nu+1}}\subset ...\end{equation}

Recall from Remark \ref{salam} that
$J^n(\bG_m)=Spf\ R[x,x^{-1},x',...,x^{(n)}]\h$ and that
$[p^{\nu}]_{J^n(\bG_m)}:J^n(\bG_m)\ra J^n(\bG_m)$
 is given at the level of rings by
$x\mapsto x^{p^{\nu}}$, $x'\mapsto \d(x^{p^{\nu}})$, etc.
 By the commutation of $J^n$ with fiber products it follows that
 $$J^n(\mu_{p^{\nu}})=Ker(J^n([p^{\nu}]_{\bG_m})=Ker([p^{\nu}]_{J^n(\bG_m)})=:J^n(\bG_m)[p^{\nu}].$$
 More generally, if $a\in R^{\times}$, and if we still denote by $a:Spec\ R \ra \bG_m$ the point defined by $x \mapsto a$ then
 $$J^n(\mu_{p^{\nu}}^a)=J^n([p^{\nu}]_{\bG_m}^{-1}(a))=(J^n([p^{\nu}]_{\bG_m})^{-1}(J^n(a))=([p^{\nu}]_{J^n(\bG_m)})^{-1}(J^n(a))$$
 where $J^n(a):Spec\ R \ra J^n(\bG_m)$
 is given, at the level of rings, by
$$x\mapsto a,\ x'\mapsto \d a,\ ...,x^{(n)}\mapsto \d^n a.$$
It follows that

\begin{proposition}
\label{plo}
We have:
$$\cO^n(\mu^a_{p^{\nu}})=\frac{R[x,x',...,x^{(n)}]\h}{(x^{p^{\nu}}-a,\d(x^{p^{\nu}})-\d a,...,\d^n(x^{p^{\nu}})-\d^n a)}.$$
In particular
$$\cO^n(\mu_{p^{\nu}})=\frac{R[x,x',...,x^{(n)}]\h}{(x^{p^{\nu}}-1,\d(x^{p^{\nu}}),...,\d^n(x^{p^{\nu}}))}.$$
\end{proposition}

Alternatively Proposition \ref{plo} follows from Lemma \ref{binomial}.

\bigskip

Let $J^n(\mu_{p^{\nu}})_1$ be the kernel of the projection $J^n(\mu_{p^{\nu}})\ra J^0(\mu_{p^{\nu}})=\widehat{\mu_{p^{\nu}}}$ and write $\cO^n(\mu_{p^{\nu}})_1:=\cO(J^n(\mu_{p^{\nu}})_1)$. Then we have the following splitting result:

\begin{proposition}
\label{smow} We have an isomorphism
$$\cO^n(\mu_p) \simeq \cO(\mu_p) \widehat{\otimes} \cO^n(\mu_p)_1$$
compatible with the group laws.
Equivalently we have
$$J^n(\mu_p)=\widehat{\mu_p} \times J^n(\mu_p)_1$$
as groups in the category of formal $p$-adic schemes over $R$.
\end{proposition}

{\it Proof}.
Let $\d^i(x^p)=F_i(x^p,x',...,x^{(i)})$; cf. Lemma \ref{45}. Then
$$\begin{array}{rcl}
\cO^n(\mu_p)_1 & = &  R[x,x',...,x^{(n)}]\h/(x-1,F_1(x^p,x'),...,F_n(x^p,x',...,x^{(n)}))\\
\  & \  & \  \\
\  & = &   R[x',...,x^{(n)}]\h/(F_1(1,x'),...,F_n(1,x',...,x^{(n)}))\end{array}$$
Hence
$$\begin{array}{rcl}
\cO^n(\mu_p) & = &  R[x,x',...,x^{(n)}]\h/(x^p-1,F_1(x^p,x'),...,F_n(x^p,x',...,x^{(n)}))\\
\  & \  & \  \\
\  & = &   R[x,x',...,x^{(n)}]\h/(x^p-1,F_1(1,x'),...,F_n(1,x',...,x^{(n)}))\\
\  & \  & \  \\
\  & \simeq  & (R[x]/(x^p-1)) \widehat{\otimes} (R[x',...,x^{(n)}]\h/(F_1(1,x'),...,F_n(1,x',...,x^{(n)})))\\
\  & \  & \  \\
\  & = & \cO(\mu_p) \widehat{\otimes} \cO^n(\mu_p)_1. \end{array}$$
The fact that the above isomorphism is compatible with group laws follows from the fact that, in view of Lemma \ref{45}, we have
$$\d^n(x_1x_2) \equiv \d^n(x_1x_2)_{|x_1=x_2=1}\ \ mod\ \ (x_1^p-1, x_2^p-1).$$
\qed

\begin{remark}
It would be interesting to understand what happens with Proposition
\ref{smow} if we considers $\mu_{p^{\nu}}$ instead of $\mu_p$.
\end{remark}

\begin{proposition}
For all $n \geq 1$ we have
$$\lim_{\stackrel{\longleftarrow}{\nu}}\cO^n(\mu_{p^{\nu}})_1=R[x',...,x^{(n)}]\h.$$
\end{proposition}

{\it Proof}.
By Proposition \ref{plo}
$$\cO^n(\mu_{p^{\nu}})_1=\frac{R[x',...,x^{(n)}]\h}{(\d(x^{p^{\nu}})_{|x=1},...,\d(x^{p^{\nu}})_{|x=1})}.$$
Since, by Lemma \ref{coo} the denominator in the last equation is in
$(p^{\nu-n+1})$ we are done by the following general well known
Lemma. \qed

\begin{lemma} If $A$ is a Noetherian ring, $I$ is an ideal,
  $A$ is $I$-adically complete,
and $(L_n)$ is a descending sequence of ideals such that $L_n
\subset I^n$ then $A=\lim_{\leftarrow} A/L_n$.\end{lemma}

{\it Proof}.  We recall the argument.
 The map $A\ra \lim_{\leftarrow}A/L_n$ is clearly injective.
 It is surjective because if
$f_n \in A$, $f_{n+1}-f_n\in L_n$ then $f_{n+1}-f_n\in I^n$ hence there exists
 $f\in A$ such that $f-f_n\in I^n$. Now fix $m$; since for $n \geq m$,
 $f-f_m=(f-f_n)+(f_n-f_m)\in I^n+L_m$ and (by \cite{mats}, Theorems 8.2
 and 8.14) $\cap_{n\geq 1} (I^n+L_m)=L_m$ we get $f-f_m\in L_m$ \qed.

\bigskip

According to our general notation (Equations \ref{genneral} and \ref{gennerall}) we next recall the rings $$\overline{\cO^n(\mu^a_{p^{\nu}})},\ \
\widetilde{\cO^n(\mu^a_{p^{\nu}})}.$$
Also recall we set $U_m=1+p^mR=U_1^{p^{m-1}}$, $m \geq 1$.

\begin{theorem}
\label{str}
Let $n,\nu \geq 1$ and $a \in U_1$.

1) If $a \not\in   U_{\nu+1}$ then
$\widetilde{\cO^n(\mu^a_{p^{\nu}})}=0$.

 2) If $a\in U_{\nu+1}$ then
\begin{equation}
\label{e1}
\widetilde{\cO^n(\mu^a_{p^{\nu}})} \simeq \widetilde{\cO^n(\mu_{p^{\nu}})}  \simeq  \frac{k[x,x',x'',...,x^{(n)}]}{((x-1)^{p^{\nu}},(x')^{p^{\nu}},...,(x^{(n)})^{p^{\nu}})}.\end{equation}
\begin{equation}
\label{e2} \overline{\cO^n(\mu^a_{p^{\nu}})}\simeq \overline{\cO^n(\mu_{p^{\nu}})} \simeq
 \frac{k[x,x',x'',...,x^{(n)}]}{(x^{p^{\nu}}-1)} \ \  \text{if $n \leq \nu$},\end{equation}
 \begin{equation}
 \label{e3}
 \overline{\cO^n(\mu^a_{p^{\nu}})}\simeq \overline{\cO^n(\mu_{p^{\nu}})}\simeq
 \frac{k[x,x',x'',...,x^{(n)}]}{(x^{p^{\nu}}-1,(x')^{p^{\nu}},...,(x^{(n-\nu)})^{p^{\nu}})} \ \  \text{if
 $n\geq \nu+1$}.\end{equation}
\end{theorem}

{\it Proof}.
 To prove assertion 1) let $a\in U_m\backslash U_{m+1}$, $1 \leq m \leq \nu$. Then $\d^m a\in R^{\times}$ hence, by Corollary \ref{modpp}, the reduction mod $p$ of $\d^m(x^{p^{\nu}})-\d^m a$ is in $k^{\times}$. By Proposition \ref{plo} $\overline{\cO^N(\mu^a_{p^{\nu}})}=0$ for all $N\geq m$ and hence $\widetilde{\cO^n(\mu^a_{p^{\nu}})}=0$ for all $n \geq 1$.

To prove assertion 2) note that the equalities (\ref{e2}) and  (\ref{e3})  follow from
 Corollary \ref{modpp} and Proposition \ref{plo}; in particular  we have
$$\lim_{\stackrel{\longrightarrow}{m}}\overline{\cO^m
(\mu_{p^{\nu}})}=\frac{k[x,x',x'',...]}{(x^{p^{\nu}}-1, (x')^{p^{\nu}}, (x'')^{p^{\nu}}, ...)}.$$
Now (\ref{e1}) follows from the fact that the intersection
$$k[x,x',...,x^{(n)}]\cap (x^{p^{\nu}}-1,(x')^{p^{\nu}}, (x'')^{p^{\nu}}, ...)$$
in the ring $k[x,x',x'',...]$ equals the ideal
$$(x^{p^{\nu}}-1,(x')^{p^{\nu}},  ...,(x^{(n)})^{p^{\nu}}).$$
\qed

\begin{corollary}
Let $n,\nu_0 \geq 1$, $a\in U_1$. Then
\begin{equation}
\label{thering}
\lim_{\stackrel{\longleftarrow}{\nu}}\widetilde{\cO^{n}(\mu^{a^{p^{\nu-\nu_0}}}_{p^{\nu}})}=
\begin{cases}
0 & \text{if $a \not\in U_{\nu_0+1}$}\\
k[[x-1,x',...,x^{(n-1)}]] & \text{if $a \in U_{\nu_0+1}$}\end{cases}\end{equation}
\end{corollary}

We also remark the following:

\begin{proposition}
\label{abv}
Assume $n, \nu \geq 1$. Then

1) $\cO^n(\mu_{p^{\nu}})$ is not a finite $R$-algebra.

2) $\cO^n(\mu_{p^{\nu}})$ is not a flat $R$-algebra.
\end{proposition}

{\it Proof}.
If $\cO^n(\mu_{p^{\nu}})$ is  a finite $R$-algebra  then $\overline{\cO^n(\mu_{p^{\nu}})}$ is  a finite $k$-algebra which is not the case, cf. equations (\ref{e2}) and (\ref{e3}). If $\cO^n(\mu_{p^{\nu}})$ is  a flat $R$-algebra then it is torsion free. Since, by equation (\ref{deltafp}), $\d(x^{p^{\nu}})\in p^{\nu}x^{p(p^{\nu}-1)}x'+p^{\nu+1}A^1$
it follows that an element in $x'+pA^1$ is zero in $\cO^n(\mu_{p^{\nu}})$ hence $x'$ is zero in
$\overline{\cO^n(\mu_{p^{\nu}})}$ hence in $\widetilde{\cO^n(\mu_{p^{\nu}})}$. Hence $x'$ is in the denominator of the ring in equation (\ref{e1}), a contradiction.
\qed

\begin{remark}
The author is indebted to A. Saha for pointing out 2) in  Proposition \ref{abv} above.
\end{remark}

\section{$p$-jets of the irreducible components of $\mu_{p^{\nu}}$}

Next note that $\mu_{p^{\nu}}$ is connected and has $\nu+1$ irreducible components:
$$\mu_{p^{\nu}}=\bigcup_{i=0}^{\nu} \mu_{p^{\nu},i},\ \ \mu_{p^{\nu},i}:=Spec\ R[\zeta_{p^i}],$$
where $\zeta_{p^i}$ is a primitive $p$-root of unity. So $\zeta_1=1$, $R[\zeta_1]=R=R[x]/(x-1)$, and
$$R[\zeta_{p^i}]=R[x]/(\Phi_{p^i}(x)),\ \ \Phi_{p^i}(x):=\frac{x^{p^i}-1}{x^{p^{i-1}}-1},\ \ i \geq 2.$$
 The scheme theoretic intersection of these  components is
$$Spec\ R[x]/(x-1,p)=Spec\ k.$$
In deep contrast with Proposition \ref{str} the $p$-jets of these components are completely uninteresting:

\begin{proposition}
\

1) $\cO^n(\mu_{p^{\nu},0})=R$ for $n \geq 1$;

2)  $\cO^n(\mu_{p^{\nu},i})=0$ for $i \geq 1$ and $n \geq 1$.
\end{proposition}

{\it Proof}.
The first equality is clear. The second follows from the fact that $\Phi_{p^i}(x+1)$ is an Eisenstein polynomial plus the following  general Proposition \ref{bel} below.
\qed

\begin{proposition}
\label{bel}
Let $f(x)=x^e+a_1x^{e-1}+...+a_{e-1}x+a_e \in R[x]$ be an Eisentein polynomial (i.e. $e \geq 2$, $a_1,...,a_e \in pR$, $a_e \not\in p^2R$) and let $X=Spec\ R[x]/(f(x))$. Then $\cO^n(X)=0$ for $n \geq 1$.
\end{proposition}

{\it Proof}.
It is enough to show that $\d f(x)$ is invertible in $\cO^1(X)\otimes k$. Now we have:
$$\begin{array}{rcl}
\d f (x) & = & \d(x^e)+\d(a_1 x^{e-1})+...+\d(a_{e-1}x)+\d a_e+C_p(x^e,a_1x^{e-1},...,a_e)\\
\  & \  & \  \\
\  & \equiv & ex^{p(e-1)}x'+ x^{p(e-1)}(\d a_1)+...+x^p(\d a_{e-1})+\d a_e\ \ mod\ \ p.
\end{array}$$
Since the image of $x$ in $\cO^1(X)\otimes k$ is  nilpotent and the image of $\d a_e$ in the same ring
is invertible it follows that the image of $\d f(x)$ in this ring is invertible which ends the proof.
\qed

\section{$p$-jets of $E[p^{\nu}]$ for ordinary elliptic curves}

We start with a review of extensions of $p^{-\nu}\bZ/\bZ$ by $\mu_{p^{\nu}}$.
For any group (respectively group scheme) $G$ we denote by $G[N]$ the kernel of the multiplication by $N$ map.
For a finite group $\Gamma$ we continue to denote by $\Gamma$  the \'{e}tale group scheme over $R$ attached to $\Gamma$; so for any connected $R$-algebra $S$, $\Gamma(S)=\Gamma$. In particular we have the connected $R$-group scheme $\mu_{p^{\nu}}=\bG_m[p^{\nu}]$. Also one can consider the \'{e}tale $R$-group scheme
$p^{-\nu}\bZ/\bZ$.
Let $R_m=R/p^mR$, $m \geq 1$. We also view, when appropriate, $\mu_{p^{\nu}}$ and $p^{-\nu}\bZ/\bZ$ as $R_m$-group schemes via base change. Then, by ``Kummer theory",
\begin{equation}\label{isooo}\begin{array}{rcl}
Ext^1_{R_m}(p^{-\nu}\bZ/\bZ,\mu_{p^{\nu}}) & \simeq & R_m^{\times}/(R_m^{\times})^{p^{\nu}}\\
\  & \simeq & (1+pR_m)/(1+pR_m)^{p^{\nu}}\\
\  & \simeq & (1+pR_m)/(1+p^{\nu+1}R_m).\end{array}\end{equation}
We will need the following explicit description of the above isomorphism.
Let $q \in 1+pR$.  Consider the finite flat $R$-scheme
$$\Gamma_{p^{\nu}}^q=\coprod_{i=0}^{p^{\nu}-1} \mu_{p^{\nu}}^{q^i}.$$
This is a group scheme with multiplication given by
$$\mu_{p^{\nu}}^{q^i}\times \mu_{p^{\nu}}^{q^j}\stackrel{can}{\longrightarrow}
\mu_{p^{\nu}}^{q^{i+j}} \stackrel{can}{\longrightarrow} \mu_{p^{\nu}}^{q^l},$$
where $0\leq l<p^{\nu}$, $i+j\equiv l$ mod $p^{\nu}$.
The functor of points of $\Gamma^q_{p^{\nu}}$ is given by
$$\Gamma^q_{p^{\nu}}(S)=\{(s,i);s\in S^{\times}, 0 \leq i <p^{\nu}, s^{p^{\nu}}=q^i\}$$
for any $R$-algebra $S$; the multiplication on points is given by $(s,i)\cdot (t,j)=(st,i+j)$ if
$i+j<p^{\nu}$ and $(s,i)\cdot (t,j)=(st/q,i+j-p^{\nu})$ if
$i+j\geq p^{\nu}$.
We have an extension
\begin{equation}
\label{theext}
0 \ra \mu_{p^{\nu}} \ra \Gamma_{p^{\nu}}^q \ra p^{-\nu}\bZ/\bZ \ra 0\end{equation}
(the second map being given on points by $(s,i)\mapsto \frac{i}{p^{\nu}}+\bZ$).
The following is ``well known":

\begin{lemma}
\label{unch}
The isomorphism (\ref{isooo}) is given by attaching to
the class of $q\in 1+pR$ in $(1+pR_m)/(1+p^{\nu+1}R_m)$ the class of the extension (\ref{theext})
in $Ext^1_{R_m}(p^{-\nu}\bZ/\bZ,\mu_{p^{\nu}})$.
\end{lemma}

Note that the system $(\Gamma_{p^{\nu}}^q;\nu\geq 1)$ is a $p$-divisible group via the morphisms
$\Gamma_{p^{\nu}}^q \ra  \Gamma_{p^{\nu+1}}^q$ given on points by $(s,i)\mapsto (s,pi)$ and given on schemes by the inclusions
$\mu_{p^{\nu}}^{q^i} \subset \mu_{p^{\nu+1}}^{q^{pi}}$.
The $p$-divisible group  $(\Gamma_{p^{\nu}}^q;\nu\geq 1)$ is an extension of the $p$-divisible group
$(p^{-\nu}\bZ/\bZ; \nu\geq 1)$ by the $p$-divisible group
  $(\mu_{p^{\nu}}; \nu\geq 1)$, where the latter are viewed as $p$-divisible groups with respect to the natural inclusions.

\bigskip

 Next we consider an elliptic curve $E/R$. References for this are \cite{sil,Katz}. Let $\overline{E}/k$ be its reduction mod $p$ and let $E^{for}$ be the formal group attached to $E$. (We use the superscript {\it for} rather than $\widehat{\ }$ because the latter is used in the present paper to denote $p$-adic completion). Let $E^{for}[p^{\nu}]$ be the kernel of the multiplication by $p^{\nu}$ on $E^{for}$, viewed as a finite flat group scheme over $R$.
Assume in what follows that $\overline{E}$ is ordinary. Then
\begin{equation}
\label{isso1}
E^{for}\simeq \bG_m^{for};\end{equation}
 we fix such an isomorphism. So we have induced isomorphisms $E^{for}[p^{\nu}]\simeq \mu_{p^{\nu}}$. Moreover we fix isomorphisms
\begin{equation}
\label{isso2}
\overline{E}(k)[p^{\nu}]\simeq \bZ/p^{\nu}\bZ\simeq p^{-\nu}\bZ/\bZ.\end{equation}
 With the isomorphisms (\ref{isso1}) and (\ref{isso2}) fixed one defines the Serre-Tate parameter $q=q(E)\in 1+pR$ of $E$ as follows.
The isomorphisms
 (\ref{isso2}) define
 a basis $(\alpha_{\nu})$ of the Tate module
$T_p\overline{E}=\lim_{\leftarrow} \overline{E}(k)[p^{\nu}]$, $\alpha_{\nu}\in \overline{E}(k)[p^{\nu}]$  a generator, $p\alpha_{\nu}=\alpha_{\nu-1}$. If $A_{\nu} \in E(R)$ lifts $\alpha_{\nu}$ then one defines the Serre-Tate parameter $q(E)\in 1+pR$ as the image of $\lim p^{\nu}A_{\nu} \in E^{for}(R)$ via the isomorphism $E^{for}(R) \simeq 1+pR$ induced by (\ref{isso1}); cf. \cite{Katz}, section 2.
On the other hand
with the isomorphisms (\ref{isso1}) and (\ref{isso2}) fixed  there are induced exact sequences of finite flat group schemes over $R$:
\begin{equation}
\label{fome}
0 \ra \mu_{p^{\nu}} \ra E[p^{\nu}]\ra p^{-\nu}\bZ/\bZ \ra 0.\end{equation}
Cf., say,  \cite{Katz}. Also by loc.cit. we have

\begin{lemma}
\label{fic}
The class of the extension (\ref{fome})  is the image of the Serre-Tate parameter $q(E)\in 1+pR$ under the isomorphism (\ref{isooo}).\end{lemma}

We conclude by Lemmas \ref{unch} and \ref{fic} that if $q=q(E)$ then $E[p^{\nu}]$ and $\Gamma^q_{p^{\nu}}$ are isomorphic as extensions over $R_m$ for any $m$; the isomorphisms are compatible as $m$ varies so we get the following:

\begin{corollary}
\label{mainnn}
$E[p^{\nu}]$ and $\Gamma^q_{p^{\nu}}$ are isomorphic as extensions over $R$ for $q=q(E)$.
  In particular if $0\leq i < p^{\nu}$ and $\theta=\frac{i}{p^{\nu}}+\bZ\in p^{-\nu}\bZ/\bZ$ then the connected component $E[p^{\nu}]_{\theta}$ of $E[p^{\nu}]$ lying above $\theta$ is isomorphic to $\mu_{p^{\nu}}^{q^i}$.
 Consequently if we fix $\nu_0$ and an index $0\leq i_0 <p^{\nu_0}$ and if $\theta=\frac{i_0}{p^{\nu_0}}+\bZ\in p^{-\nu_0}\bZ/\bZ$ then the inductive system
$$E[p^{\nu_0}]_{\theta}\subset E[p^{\nu_0+1}]_{\theta} \subset E[p^{\nu_0+2}]_{\theta} \subset ...\subset
E[p^{\nu}]_{\theta}\subset ...$$
identifies with the inductive system
$$\mu_{p^{\nu_0}}^{q^{i_0}} \subset \mu_{p^{\nu_0+1}}^{(q^{i_0})^p}\subset \mu_{p^{\nu_0+2}}^{(q^{i_0})^{p^2}}
\subset ... \subset \mu_{p^{\nu}}^{(q^{i_0})^{p^{\nu-\nu_0}}}\subset ...$$\end{corollary}

Putting together Theorem \ref{str} and Corollary \ref{mainnn} (and making the change of variables $x \mapsto x+1$) we get:

\begin{theorem}
\label{painnn}
Let $E/R$ be an elliptic curve with ordinary reduction and Serre-Tate parameter $q=q(E)\in U_1$.
Let $n, \nu_0 \geq 1$,  $\theta\in p^{-\nu_0}\bZ/\bZ$, and assume $\nu_0$ is minimal with this property.
Let $E[p^{\nu}]_{\theta}$ be the connected component of $E[p^{\nu}]$ lying over $\theta$. Then:

1) If $q\not\in U_{\nu_0+1}$ and $\theta \neq 0$ then
$\widetilde{\cO^n(E[p^{\nu}]_{\theta})}=0$.

 2) If $q\in U_{\nu_0+1}$ or $\theta = 0$ then
\begin{equation}
\label{e1}
\widetilde{\cO^n(E[p^{\nu}]_{\theta})}   \simeq  \frac{k[x,x',x'',...,x^{(n)}]}{(x^{p^{\nu}},(x')^{p^{\nu}},...,(x^{(n)})^{p^{\nu}})}.\end{equation}
\begin{equation}
\label{e2} \overline{\cO^n(E[p^{\nu}]_{\theta})} \simeq
 \frac{k[x,x',x'',...,x^{(n)}]}{(x^{p^{\nu}})} \ \  \text{if $n \leq \nu$},\end{equation}
 \begin{equation}
 \label{e3}
 \overline{\cO^n(E[p^{\nu}]_{\theta})}\simeq
 \frac{k[x,x',x'',...,x^{(n)}]}{(x^{p^{\nu}},(x')^{p^{\nu}},...,(x^{(n-\nu)})^{p^{\nu}})} \ \  \text{if
 $n\geq \nu+1$}.\end{equation}
\end{theorem}

\begin{corollary}
Let $n,\nu_0 \geq 1$, $a\in U_1$. Then
\begin{equation}
\label{thering}
\lim_{\stackrel{\longleftarrow}{\nu}}\widetilde{\cO^{n}(E[p^{\nu}]_{\theta})}=
\begin{cases}
0 & \text{if $q\not\in U_{\nu_0+1}$ and $\theta \neq 0$ }\\
k[[x,x',...,x^{(n-1)}]] & \text{if $q\in U_{\nu_0+1}$ or $\theta = 0$. }\end{cases}\end{equation}
\end{corollary}

Also Corollaries \ref{abv} and \ref{mainnn} imply

\begin{corollary}
\label{kamel} Let $n,\nu \geq 1$. Then  the map
$$J^n([p^{\nu}]_E)=[p^{\nu}]_{J^n(E)}:J^n(E)\ra J^n(E)$$ is neither finite not flat.
\end{corollary}

An analogue of Corollary \ref{aftert} should hold nevertheless in
the elliptic case as well.

\section{$p$-jets of $\cF[p^{\nu}]$}

Recall that a series $F(x)\in xR[[x]]$ without constant term is said to have finite height if $F \not\equiv 0$ mod $p$; if this is the case the height of $F$ is defined as the largest integer $h \geq 0$ such that $F \in R[[x^{p^h}]]+pR[[x]]$.

\begin{remark}
The main example we have in mind here arises as follows.
Consider a formal group law $\cF \in R[[x_1,x_2]]$. By \cite{sil} the multiplication by $p$ in $\cF$ is given by a series
$F(x):=[p]_{\cF}(x)$ satisfying $F(x)\equiv px$ mod $x^2$. The height of $\cF$ is defined to be the height of $F(x)$; if the height is finite then it is $\geq 1$. Not every series of height $\geq 1$ which is $\equiv px$ mod $x^2$ is the multiplication by $p$ of a formal group law $\cF$; indeed if $F(x)=[p]_{\cF}$ has height $h$ then one knows that the $x$-adic valuation of the reduction mod $p$ of $F$ in $k[[x]]$ is exactly $p^h$; cf. \cite{sil}, p.127.
\end{remark}

Let $F(x)\in xR[[x]]$ be a series of finite height $h\geq 1$ and let $F^{\circ \nu}=F \circ ... \circ F$ be the $\nu$ fold composition of $F$ with itself for $\nu\geq 1$. Let $ep^h$ be the $x$-adic valuation of the reduction mod $p$ of $F$. (So if $F(x)=[p]_{\cF}(x)$ for some formal group law $\cF$ then $e=1$.)
By ``Weierstrass preparation" (cf. \cite{Lang}, p. 130)   $F^{\circ \nu}=U_{\nu}\cdot P_{\nu}$ where $U_{\nu}\in R[[x]]^{\times}$ and $P_{\nu}\in R[x]$ is monic of degree $e^{\nu}p^{h\nu}$, $P_{\nu}\equiv x^{e^{\nu}p^{h\nu}}$ mod $p$.
Consider the scheme:
$$X_{\nu}:=Spec\ \frac{R[[x]]}{(F^{\circ \nu})}=Spec\ \frac{R[[x]]}{(P_{\nu})}\simeq Spec\ \frac{R[x]}{(P_{\nu})};$$
the latter isomorphism follows from ``Euclid division"  by $P_{\nu}$ in $R[[x]]$; cf. \cite{Lang}, p. 129.
So $X_{\nu}$ is a finite flat scheme over $R$ of degree $e^{\nu}p^{h\nu}$ and we have a natural sequence of closed immersions
\begin{equation}
\label{ciine}
X_1\subset X_2 \subset ... \subset X_{\nu} \subset ...  \end{equation}
 Note that if $F(x)=[p]_{\cF}(x)$ is the multiplication by $p$ on  some formal group law $\cF$ then $X_{\nu}=\cF[p^{\nu}]$, where the latter is the kernel of $[p^{\nu}]_{\cF}$ on $\cF$ and indeed the  inductive system (\ref{ciine}) coincides with the $p$-divisible group
 $$\cF[p]\subset \cF[p^2] \subset ... \subset \cF[p^{\nu}]\subset ...$$
  of $\cF$; cf. \cite{Tate}.
According to  our general notation in equations \ref{genneral} and \ref{gennerall} we consider the rings
\begin{equation}
\label{alco}
\overline{\cO^n(X_{\nu})}, \ \widetilde{\cO^n(X_{\nu})}.
\end{equation}

\begin{remark}
Recall (cf. \cite{Lubin}, p. 480) that any formal group law $\cF$ over $R$ of height $h=1$ is isomorphic to the multiplicative formal group law hence in particular $\cF[p^{\nu}]\simeq \mu_{p^{\nu}}$.
Hence Theorem \ref{str} provides a computation of the rings (\ref{alco}) in this case. The following Theorem covers the case of formal group laws $\cF$ of height $h \geq 2$. More generally the Theorem below treats the case of iterates of series of height $\geq 1$ which are not necessarily coming from formal groups; so even for height $1$ the Theorem below is not covered by Theorem \ref{str}.
\end{remark}

\begin{theorem}
\label{mor}
Let $F(x)\in xR[[x]]$ be a series of finite height $h \geq 1$ satisfying $F(x)\equiv px$ mod $x^2$.
 For all $\nu\geq 1$ consider the scheme $X_{\nu}:=Spec\ \frac{R[[x]]}{(F^{\circ \nu})}$.
 Then we have:
\begin{equation}
\label{secon}
\widetilde{\cO^n(X_{\nu})}=
\frac{k[x,x',...,x^{(n)}]}{(x^{p^{\nu}},(x')^{p^{\nu}},...,(x^{(n)})^{p^{\nu}})}\ \ \ \text{if $n\geq 1$},
\end{equation}
\begin{equation}
\label{firs}
\overline{\cO^n(X_{\nu})}=\frac{k[x,x',...,x^{(n)}]}{(x^{p^{\nu}},
(x')^{p^{\nu}},...,(x^{(n-\nu)})^{p^{\nu}})}\ \ \ \text{if $n\geq \nu$},
\end{equation}
\begin{equation}
\overline{\cO^n(X_{\nu})}=\frac{k[x,x',...,x^{(n)}]}{(x^{p^{\mu}})}\ \ \ \text{if $1 \leq n \leq \nu-1$},
\end{equation}
where $\mu\geq \nu$.
\end{theorem}

\begin{remark}
Let $\cF$ be a formal group law over $R$ of finite height $h$. Note that the morphism
$$\overline{\cO(\cF[p^{\nu}])}\ra \overline{\cO^n(\cF[p^{\nu}])}$$
is  injective if $h=1$ (cf. Theorem \ref{str}) but not injective if the height $h \geq 2$ and $n \geq \nu$ (cf. Theorem \ref{mor}).
\end{remark}

\begin{remark}
To prove the theorem above note first that:
\begin{equation}
\label{cd}
\begin{array}{rcl}
\cO^n(X_{\nu}) & = & R[x,x',...,x^{(n)}]\h/(P_{\nu},\d P_{\nu},...,\d^n P_{\nu})\\
\  & \  & \  \\
\  & = & R[[x]][x',...,x^{(n)}]\h/(P_{\nu},\d P_{\nu},...,\d^n P_{\nu})\\
\ & \  & \  \\
\  & = & R[[x]][x',...,x^{(n)}]\h/(F^{\circ \nu},\d (F^{\circ \nu}),...,\d^n (F^{\circ \nu})).\end{array}\end{equation}
So we will concerned from now on with understanding the structure of the expressions $\d^i(F^{\circ \nu})$. To do this we need to develop some filtration machinery on power series.
\end{remark}

We start by  considering  the decreasing filtration of $\cA^0:=R[[x]]$ by the subrings $\cA^0_{\nu}$, $\nu\geq 1$, defined by
$$\cA^0_{\nu}=R[[x^{p^{\nu}}]]+pR[[x^{p^{\nu-1}}]]+p^2R[[x^{p^{\nu-2}}]]+...+p^{\nu}R[[x]]\subset R[[x]].$$
Let $v_p$ be the $p$-adic valuation on $R$.

\begin{lemma}
\label{hess}

\

1) $\cA^0_{\nu}=\{\sum_{n \geq 0} a_n x^n\in R[[x]]\ ;\ v_p(a_n)\geq \nu-v_p(n)\}$.

2) If $G_1,G_2,G_3,...\in \cA^0_{\nu}$, $G_m\in x^m R[[x]]$. Then $\sum_{m \geq 1} G_m \in \cA^0_{\nu}$.

3) If $H \in \cA^0_{\nu}$, $H(0)=0$, and $G \in R[[x]]$ then $G(H(x)) \in \cA^0_{\nu}$.

4) $p\cA^0_{\nu}\subset \cA^0_{\nu+1}$.

5) If $G\in \cA^0_{\nu}$ then $G^{p}\in \cA^0_{\nu+1}$.

6) If $F\in \cA^0_1$ and  $F(0)=0$ then $F^{\circ \nu}\in \cA^0_{\nu}$.
\end{lemma}

{\it Proof}. 1) is easy.
2) clearly follows from 1).
3) clearly follows from 2).
4) and 5) are clear. 6) follows from 3), 4), 5).
\qed

\bigskip

We continue by considering the filtration
$$\cA^n=R[[x]][x',...,x^{(n)}]\h,\ \ n \geq 0$$
on
$$\cA:=\bigcup_{n\geq 0} \cA^n.$$
(Here $\cA^0=R[[x]]$.) There is a natural $p$-derivation $\d$ on $\cA$ sending $\d x=x'$, $\d x'=x''$, etc. Note that $\d \cA^{n}\subset \cA^{n+1}$ for all $n$. So according to
equation (\ref{AAA}) we may then consider the filtration $\cA^{\{n\}}$ on $\cA$. Finally let
$$\cI=(x,x',x'',...)\subset \cA.$$
So we may consider the descending filtration of $\cI$ by ideals $\cI^{[p^{\nu}]}$, $\nu\geq 0$.
Note that with $A^n,A,I$ as in \ref{oz} we have $A^n\subset \cA^n$, $A \subset \cA$, $I\subset \cI$, and hence $I^{[p^{\nu}]}\subset \cI^{[p^{\nu}]}$.

\begin{lemma}
\label{lem4}
Let $n,i\geq 1$.
 $$\d^n(p^ix)\in
 \begin{cases}
 p^{i-n} \phi^n(x)+(p^{i-n+1}\cA^{\{0\}})\cap \cI^{[p^{i+1}]} & \text{if $n\leq i$}\\
 \phi^i(x^{(n-i)})+\cA^{\{n-i-1\}}\cap\cI^{[p^{i+1}]} & \text{if $n \geq i+1$}\end{cases}$$
\end{lemma}

{\it Proof}.
Induction on $n$. The case $n=1$ is clear. Now assume the above is true for some $n \geq 1$.
If $n\leq i-1$ then by Lemma \ref{lem2} $p^{i-n}\phi(x)^p\in p^{i-n}\cI^{[p^{n+1}]}\subset \cI^{[p^{i+1}]}$ so
$$\begin{array}{rcl}
\d^{n+1}(p^ix) & \in & \d(p^{i-n}\phi^n(x))+ \d((p^{i-n+1}\cA^{\{0\}})\cap \cI^{[p^{i+1}]})\\
\  & \  & +C_p(p^{i-n}\cA^{\{0\}},(p^{i-n+1}\cA^{\{0\}})\cap \cI^{[p^{i+1}]})\\
\  & \subset &  p^{i-n-1}\phi^{n+1}(x) -p^{(i-n)p-1}\phi^n(x)^p+(p^{i-n}\cA^{\{0\}})\cap \cI^{[p^{i+1}]}\\
\  & \  & + (p^{i-n+1}\cA^{\{0\}})\cap \cI^{[p^{i+1}]}\\
\  & \subset & p^{i-n-1}\phi^{n+1}(x)
+(p^{i-n}\cA^{\{0\}})\cap \cI^{[p^{i+1}]}.\end{array}$$
If $n \geq i+1$ we have
$$\begin{array}{rcl}\d^{n+1}(p^ix) & \in & \d(\phi^i(x^{(n-i)}))+ \d(\cA^{\{n-i-1\}})\cap \cI^{[p^{i+1}]})\\
\  & \  & +C_p(\cA^{\{n-i\}},\cA^{\{n-i-1\}}\cap \cI^{[p^{i+1}]})\\
\  & \in & \phi^i(x^{(n+1-i)}) +\cA^{\{n-i\}}\cap \cI^{[p^{i+1}]}.\end{array}
$$
The case $n= i$ is similar.
\qed

\bigskip

\begin{lemma}
\label{lem5}
Let $G\in xR[[x]]$ be any series without constant term. Then for all $n \geq 1$:
$$\d^n(G(x))\in \left( \frac{dG}{dx} \right)^{p^n}x^{(n)}+\cA^{\{n-1\}}\cap \cI^{[p]}.$$
\end{lemma}

{\it Proof}.
We proceed by induction on $n$. To check the statement for $n=1$ write $G(x)=\sum_{n \geq 1} a_n x^n$; then
$$\begin{array}{rcl}
\d(G(x)) & = & \frac{1}{p}[ \sum_{n\geq 1} \phi(a_n)(x^p+px')^n- (\sum_{n\geq 1}a_nx^n)^p]\\
\  & \  & \  \\
\  & \subset & x^p\cA^0
+(\sum_{n\geq 1} na_nx^{n-1})^px'+px'\cA^1\\
\  & \  & \  \\
\  & \subset & \left( \frac{dG}{dx} \right)^{p}x'+\cA^{\{0\}}\cap \cI^{[p]}, \end{array}$$
which settles the case $n=1$. For the induction step, assuming the statement true for some $n\geq 1$, we have
$$\begin{array}{rcl}
\d^{n+1}(G(x)) & \in & \d\left(\left( \frac{dG}{dx} \right)^{p^n}x^{(n)}\right)
+ \d( \cA^{\{n-1\}}\cap \cI^{[p]})+C_p(\cA^{\{n\}},\cA^{\{n-1\}}\cap \cI^{[p]})\\
\  & \  & \  \\
\  & \in & \left( \frac{dG}{dx} \right)^{p^{n+1}}x^{(n+1)}+
\phi(x^{(n)}) \cdot \d \left(\left( \frac{dG}{dx} \right)^{p^{n}}\right) + \cA^{\{n\}}\cap \cI^{[p]}. \end{array}$$
Now, using Theorem \ref{coo}, we have
$$\begin{array}{rcl}
\phi(x^{(n)}) \cdot \d \left(\left( \frac{dG}{dx} \right)^{p^{n}}\right) & \in &
(\cA^n \cap \cI) \cdot p^n\cA^1\\
\  & \subset & \cA^n \cap p^n\cI\\
\  & \subset & \cA^{\{n\}} \cap \cI^{[p]}\end{array}$$
and we are done.
\qed

\bigskip

Next let $\Sigma=\Sigma(x)\in xR[[x]]$ be any series without constant term and consider the unique ring endomorphism
$\Sigma^*:\cA \ra \cA$ such that $\Sigma^* x=\Sigma(x)$, $\Sigma^* x'=\d (\Sigma(x))$, $\Sigma^* x''=\d^2 (\Sigma(x))$, etc. Clearly $\Sigma^*(\cA^n) \subset \cA^n$ for $n \geq 0$ and hence
$\Sigma^*(\cA^{\{n\}}) \subset \cA^{\{n\}}$ for $n \geq 0$.
It is trivial to see that $\Sigma^*$ and $\d$ commute on $\cA$; similarly $\Sigma^*$ and $\phi$ commute on $\cA$. Moreover for any two series $\Sigma_1,\Sigma_2$ we have the following compatibility of upper $*$ with composition: $(\Sigma_1\circ \Sigma_2)^*=\Sigma_2^* \circ \Sigma_1^*$.

For any integer $a\in \bZ$ write $a^+=\max\{a,0\}$.

\begin{lemma}
\label{lem6}
Assume $\Sigma(x)=x^{p^m}$, $m \geq 1$, $\nu\geq 0$, $n \geq 0$; then:

1) $\Sigma^* \cI^{[p^{\nu}]}\subset \cI^{[p^{\nu+m}]}$;

2) $\Sigma^* \cA^{\{n\}}\subset \cA^{\{(n-m)^+\}}$.
\end{lemma}

{\it Proof}. To check 1)
it is enough to check it for $m=1$. Now 1)
follows from the following computations in which $i+j=\nu$:
$$\Sigma^* (p^i(x^{(s)})^{p^j})=p^i(\d^s(x^p))^{p^j}\in
 p^i \cI^{[p^{j+1}]} \subset \cI^{[p^{i+j+1}]}=\cI^{[p^{\nu+1}]};$$
in the above we used the fact that since
$x^p\in \cI^{[p]}$ we have $\d^s(x^p)\in \cI^{[p]}$ hence
$(\d^s(x^p))^{p^j}\in \cI^{[p^{j+1}]}$; cf. Lemma \ref{lem1}.
To prove 2)  it is enough, by the compatibility with composition, to prove these two statements for $m=1$ and $n \geq 1$ which we now assume. Now by Theorem \ref{coo} we have $\d^n(x^p)\in \cA^{\{n-1\}}$ for $n \geq 1$. Consequently, for $n\geq 1$ and $F \in \cA^{n+i}$ we have
 $$\Sigma^*(p^i F(x,...,x^{(n+i)}))= p^i F(x^p,...,\d^{n+i}(x^p))\in p^i \cA^{\{n+i-1\}}\subset \cA^{\{n-1\}}$$
 which proves 2).
\qed

\begin{lemma}
\label{lem7}
Let $G(x)\in xR[[x]]$, $m\geq 1$, $\nu\geq 1$. Then for $i+j=\nu$, $i,j\geq 0$ we have
that $\d^m(p^iG(x^{p^j}))$ belongs to
\begin{equation}
\label{dream1}
\begin{cases}
p^{i-m} \phi^m(G(x^{p^{j}}))+(p^{i-m+1}\cA^{\{0\}})\cap \cI^{[p^{\nu+1}]},&  \text{if $m \leq i$}\\
\phi^i\{\left( \frac{dG}{dx}(x^{p^{j}})\right)^{p^{m-i}} \d^{m-i}(x^{p^{j}})\}
+\cA^{\{(m-\nu-1)^+\}} \cap \cI^{[p^{\nu+1}]},  & \text{if $m \geq i+1$}.
\end{cases}\end{equation}
\end{lemma}

{\it Proof}. Set $\Sigma(x)=x^{p^{j}}$.
Using Lemmas \ref{lem4}, \ref{lem5}, \ref{lem6} we have the following computation for $1 \leq m \leq  i$:
$$
\begin{array}{rcl}
\d^m(p^iG(x^{p^{j}})) & = & \d^m(\Sigma^* G^* (p^ix))\\
\  & = & \Sigma^* G^* (\d^m(p^ix))\\
\  & = & \Sigma^* G^*\{p^{i-m} \phi^m(x)+(p^{i-m+1}\cA^{\{0\}})\cap \cI^{[p^{i+1}]}\}\\
\  & = & \Sigma^* G^*(p^{i-m} \phi^m(x))+(p^{i-m+1}\cA^{\{0\}})\cap \cI^{[p^{i+1+j}]}\\
\  & = & p^{i-m} \phi^m(G(x^{p^{j}}))+(p^{i-m+1}\cA^{\{0\}})\cap \cI^{[p^{\nu+1}]}.
\end{array}
$$
For $m \geq i+1$  we have:
$$\begin{array}{rcl}
\d^m(p^iG(x^{p^{j}})) & = &  \Sigma^* G^* (\d^m(p^ix))\\
\  & \in & \Sigma^* G^* \{\phi^i(x^{(m-i)})+\cA^{\{m-i-1\}} \cap \cI^{[p^{i+1}]}\}\\
\  & = & \Sigma^* G^* \{\phi^i(x^{(m-i)})\}+\cA^{\{(m-i-1-j)^+\}} \cap \cI^{[p^{i+1+j}]}\\
\  & = & \phi^i\Sigma^* (\d^{m-i} G)+ \cA^{\{(m-\nu-1)^+\}} \cap \cI^{[p^{\nu+1}]}\\
\  & = & \phi^i\Sigma^* \{ \left(\frac{dG}{dx}\right)^{p^{m-i}} x^{(m-i)}+\cA^{\{m-i-1\}} \cap \cI^{[p]} \}\\
\  & \  & + \cA^{\{(m-\nu-1)^+\}} \cap \cI^{[p^{\nu+1}]}\\
\  & = & \phi^i\{\left( \frac{dG}{dx}(x^{p^{j}})\right)^{p^{m-i}}\cdot \d^{m-i}(x^{p^{j}})\}
+\cA^{\{(m-\nu-1)^+\}} \cap \cI^{[p^{\nu+1}]}.
\end{array}$$
\qed

\bigskip

For any element $f\in \cA^n$ we denote by $\overline{f}\in \overline{\cA}^n:=\overline{\cA^n} =k[[x]][x',...,x^{(n)}]$ the class of $f$ mod $p$. In particular one can consider the reduction mod $p$ of the ideal $\cI$, $\overline{\cI}=(x,x',x'',...)\subset \overline{\cA}=k[[x]][x',x'',...]$. Then
$$\overline{\cI^{[p^{\nu}]}}=\overline{\cI}^{[p^{\nu}]}=(x^{p^{\nu}}, (x')^{p^{\nu}},(x'')^{p^{\nu}},...).$$
Moreover clearly
$$\overline{\cA}^n\cap \overline{\cI}^{[p^{\nu}]}=(x^{p^{\nu}}, (x')^{p^{\nu}},...,(x^{(n)})^{p^{\nu}})\subset \overline{\cA}^n.$$
Also for any $f\in \cA^{\{n\}}\cap \cI^{[p^{\nu}]}$ we have $\overline{f} \in \overline{\cA}^n\cap\overline{\cI}^{[p^{\nu}]}$.
Now Lemma \ref{lem7} trivially implies :

\begin{lemma}
\label{lem8}
Let $m\geq 0$, $\nu\geq 1$, $i+j=\nu$, $i,j\geq 0$, $G\in xR[[x]]$. Then the element
$$\overline{\d^m(p^iG(x^{p^j}))}\in \overline{\cA}^m$$
is given by
$$\begin{cases}
0, & \text{if $m \leq i-1$}\\
G(x^{p^j})^{p^i}, & \text{if $m=i$}\\
\text{(elt. in $\overline{\cA}^0 \cap \overline{\cI}^{[p^{\nu+1}]}$)}, & \text{if $i+1\leq m\leq \nu$}\\
\left(x^{p^j-1}\frac{dG}{dx}(x^{p^j})\right)^{p^m}  (x^{(m-\nu)})^{p^{\nu}}+\text{(elt. in $\overline{\cA}^{m-\nu-1}\cap\overline{\cI}^{[p^{\nu+1}]}$)}, & \text{if $m \geq \nu+1$}.\end{cases}$$
In particular
$$\overline{\d^m(p^iG(x^{p^j}))}\in \overline{\cA}^{(m-\nu)^+}\cap \overline{\cI}^{[p^{\nu}]}.$$
\end{lemma}

\begin{lemma}
\label{lem9}
Let $y$ be an $N$-tuple $y_1,...,y_N$ of variables. Then, for $n\geq 1$,
$$\d^n\left(\sum_{i=1}^N y_i\right)=\sum_{i=1}^N y_i^{(n)} + P_{N,n}(y,y',...,y^{(n-1)})$$
in $R\{y\}$ where $P_{N,n}$ is a polynomial with $\bZ$-coefficients without constant term or linear terms.
\end{lemma}

{\it Proof}. Induction on $n$.
\qed

\bigskip

For the next Lemma note that for any series $G(x)\in \cA^0_{\nu}$ we
have $\frac{dG}{dx}\in p^{\nu}R[[x]]$ so we may consider the series
$$\overline{p^{-\nu}\frac{dG}{dx}}\in k[[x]].$$
Then Lemmas \ref{lem8} and \ref{lem9} immediately imply:

\begin{lemma}
\label{lem10} Consider a series
$G(x)=\sum_{j=0}^{\nu}p^{\nu-j}G_j(x^{p^j})$, $G_j(x)\in xR[[x]]$.
Then, for $\nu\geq 1$ and $n \geq 0$ the element $\overline{\d^n
G}\in \overline{\cA}^n$ is given by
$$
\overline{\d^n G}=\begin{cases}
(G_{\nu-n}(x^{p^{\nu-n}}))^{p^n}+\text{(elt. in $x^{2p^{\nu}}k[[x]]$)} & \text{if $0 \leq n \leq \nu$} \\
\  & \ \\
\left( \overline{p^{-\nu}\frac{dG}{dx}}\right)^{p^n}
(x^{(n-\nu)})^{p^{\nu}}+\text{(elt. in $\overline{\cA}^{n-\nu-1}\cap
\overline{\cI}^{[\nu]}$)} & \text{if $n\geq \nu+1$}.
\end{cases}
$$
\end{lemma}

\bigskip

{\it Proof of Theorem \ref{mor}}. By assertion 6) in Lemma
\ref{hess} we may write $F^{\circ
\nu}=\sum_{j=0}^{\nu}p^{\nu-j}G_j(x^{p^{j}})$,  $G_j \in R[[x]]$, $j
\geq 0$. We may choose the $G_j$s in $xR[[x]]$ and then
$G_0(x)\equiv x$ mod $x^2$. Also
$\overline{p^{-\nu}\frac{dG}{dx}}\equiv 1$ mod $x$. We conclude by
Lemma \ref{lem10} and equation (\ref{cd}).\qed

\section{$p$-jets of $E[p^{\nu}]$ for supersingular elliptic curves}

\begin{theorem}
Let $E/R$ be an elliptic curve with supersingular reduction and $E[p^{\nu}]$ the kernel of the multiplication by $p$. Then for any $n \geq 1$ we have:
\begin{equation}
\label{secon}
\widetilde{\cO^n(E[p^{\nu}])}=
\frac{k[x,x',...,x^{(n)}]}{(x^{p^{\nu}},(x')^{p^{\nu}},...,(x^{(n)})^{p^{\nu}})}\ \ \ \text{if $n\geq 1$},
\end{equation}
\begin{equation}
\label{firs}
\overline{\cO^n(E[p^{\nu}])}=\frac{k[x,x',...,x^{(n)}]}{(x^{p^{\nu}},
(x')^{p^{\nu}},...,(x^{(n-\nu)})^{p^{\nu}})}\ \ \ \text{if $n\geq \nu$},
\end{equation}
\begin{equation}
\overline{\cO^n(E[p^{\nu}])}=\frac{k[x,x',...,x^{(n)}]}{(x^{p^{\mu}})}\ \ \ \text{if $1 \leq n \leq \nu-1$},
\end{equation}
where $\mu \geq \nu$.
\end{theorem}

{\it Proof}.
Since $E$ has supersingular reduction
$E[p^{\nu}]$ is connected so it is equal to $\cF[p^{\nu}]$ where $\cF$ is the formal group law of $E$ and we conclude by Theorem \ref{mor}.
\qed

\bibliographystyle{amsplain}

\end{document}